\documentclass[a4paper,12pt,reqno]{amsart}

\usepackage{amsmath,amssymb,amsthm,latexsym,mathrsfs,eucal}

\newtheorem{thm}{Theorem}[section]
\newtheorem{assum}[thm]{Assumption}

\newtheorem{lemma}[thm]{Lemma}
\newtheorem{cor}[thm]{Corollary}

\theoremstyle{definition}
\newtheorem{defn}[thm]{Definition}

\DeclareMathOperator{\im}{Im}
\DeclareMathOperator{\End}{End}
\newcommand{\nc}{\newcommand}
\newcommand{\rnc}{\renewcommand}
\renewcommand{\l}{\lambda}
\renewcommand{\d}{\delta}
\rnc{\k}{\kappa}

\newcommand{\noi}{\noindent}
\newcommand{\bnk}{\mathscr{B}_n^k}
\newcommand{\bt}{\mathscr{B}_2^k}
\newcommand{\ix}{X^{-1}}
\newcommand{\iy}{Y^{-1}}
\rnc{\b}[1]{\,\overline{\!#1}}
\nc{\e}{\epsilon}
\nc{\kum}[1]{\sum_{#1=0}^{k-1}}
\nc{\kkum}[1]{\sum_{#1=0}^{k}}
\nc{\tw}[1]{0 \leq #1 \leq k-1}
\nc{\tf}{\text{,\quad for all }}
\nc{\iif}{\text{,\quad if }}
\nc{\la}[1]{\langle #1 \rangle}
\rnc{\=}{&=&}
\nc{\rel}[1]{\stackrel{(\ref{#1})}{=}}
\nc{\eq}[1]{\stackrel{\hphantom{(\ref{#1})}}{=}}

\newcommand{\be}{\begin{align*}}
\newcommand{\ee}{\end{align*}}
\rnc{\=}{&=&}
\newcommand{\beq}{\begin{align}}
\newcommand{\eeq}{\end{align}}

\rnc{\thefootnote}{\fnsymbol{footnote}}

\begin{document}

\begin{titlepage}
\title[The cyclotomic BMW algebra $\bt$]{The cyclotomic BMW algebra associated with the two string type $\boldsymbol{B}$ braid group}
\author{Stewart Wilcox}
\address{Stewart Wilcox \\ Department of Mathematics \\ Harvard University \\ 1 Oxford Street \\ Cambridge MA 02318 \\ USA}
\email{swilcox@fas.harvard.edu}
\author{Shona Yu}
\address{Shona Yu \\ Faculteit Wiskunde en Informatica \\ Technische Universiteit Eindhoven \\ 5600MB Eindhoven \\ The Netherlands}
\email{s.h.yu@tue.nl}
\footnotetext{Date: June, 2007. Available from \texttt{http://arxiv.org/abs/math/0611518}}
\end{titlepage}

\begin{abstract}
The cyclotomic Birman-Murakami-Wenzl (or BMW) algebras $\bnk$, introduced by R.~H\"aring-Oldenburg,
are extensions of the cyclotomic Hecke algebras of Ariki-Koike, in the same way as the
BMW algebras are extensions of the Iwahori-Hecke algebras of type~$A$.
In this paper we produce a basis of the two strand algebra $\bt$ and construct its left
regular representation. In the process we obtain relations amongst the parameters which 
play an important role in the $n$ strand case.
\end{abstract}
\maketitle

\section{Introduction}
\setcounter{page}{1}
Motivated by knot theory associated to the Artin braid group of type $A$ and the Kauffman link polynomial,  Murakami \cite{M87} and Birman and Wenzl \cite{BW89} defined what are now known as the BMW algebras.  
The Artin braid group relations of type~$A$ appear in the defining relations of the BMW algebras. Also, the Iwahori-Hecke algebra of the symmetric group (the Coxeter group of type $A$) appears naturally as a quotient of this algebra.
Motivated by type $B$ knot theory and the Ariki-Koike algebras~\cite{AK94} (also known as the cyclotomic Hecke algebras of type $G(k,1,n)$), H\"aring-Oldenburg \cite{HO01} introduced the ``cyclotomic BMW algebras" $\bnk$, of which the Ariki-Koike algebras appear naturally as a quotient in a similar fashion.

Unlike the BMW and Ariki-Koike algebras, one needs to impose extra so-called ``admissibility conditions" (see Definition \ref{defn:adm}) on 
the parameters of the ground ring in order for the algebras $\bnk$ to be well-behaved. This is because the polynomial relation of order 
$k$ imposed on $Y$ (which does not appear in the BMW or affine BMW algebras) potentially causes torsion on 
elements associated with certain tangles on two strands (which are not braid diagrams, and so do not appear in the Ariki-Koike algebras). 
These conditions may be determined precisely by focussing on the algebra $\bt$ and its representations.
Remarkably these conditions, which are required for $\bt$ to be free of the expected rank $3k^2$, ensure that $\bnk$ is free of the expected rank $k^n (2n-1)!!$
for any positive integer $n$. Presently, under the assumption that a certain element of the ground ring is not a zero divisor (see Assumption \ref{NAZD}), we describe these admissibility conditions explicitly and prove that, for a ring $R$ with parameters satisfying these conditions, the algebra $\bt(R)$ is $R$-free of rank $3k^2$. In \cite{bnk1,Yu07} we remove Assumption \ref{NAZD}, and modify the relations of Definition \ref{defn:adm} to produce relations which are \emph{necessary and sufficient} for the freeness results to hold in the general case for \emph{all} $n$, allowing us to construct a ``generic" admissible ground ring.

The structure of the paper is as follows. In Section 2, we give the definition of $\bnk$, and describe the class of ground rings we will consider. In Section 3, we give explicitly derive the admissibility conditions, which naturally arise as we construct a particular $\bt$-module $V$ of rank $k$ (see Lemma \ref{lem:V}). Then in Section 4, using this rank $k$ module, we are able to construct the regular representation of $\bt$ and provide an explicit basis of the algebra.
These results are stated but incompletely proved in H\"aring-Oldenburg~\cite{HO01}; specifically, additional arguments are needed to prove Lemma~25 of \cite{HO01}. In the present paper we take a slightly different approach and the arguments we offer correct this problem. Goodman and Hauschild Mosley \cite{GH109,GH209} and Rui and Xu \cite{Rui} have proven similar results using different methods, and imposing different conditions on the ground ring. We will show (see Corollary \ref{generic}) that the conditions we impose are the most general under which $\bt$ is free of the expected rank.

The authors thank Bob Howlett for his helpful suggestions during the preparation of this manuscript, Fred Goodman for useful email correspondences, and the referee for his/her comments on this paper.
%\newpage
\section{The Cyclotomic BMW algebras}
\begin{defn} \label{defn:bnk}
Fix natural numbers $n$ and $k$. Let $R$ be a commutative ring containing units $q_0,q,\l$ and further 
elements $A_0, A_1, \ldots, A_{k-1}, q_1, q_2,\ldots, q_{k-1}$ such that $\l - \l^{-1} = \d(1-A_0)$ holds, where $\d := q - q^{-1}$. \\ 
The \emph{\textbf{cyclotomic BMW algebra}} $\bnk(R)$ is defined to be the unital associative $R$-algebra 
generated by $Y^{\pm 1}, X_1^{\pm 1}, \ldots, X_{n-1}^{\pm 1}$ and 
$e_1, \ldots, e_{n-1}$ subject to the following relations:

\begin{eqnarray*}
X_i - X_i^{-1} &=& \d(1-e_i) \\
X_iX_j &=& X_jX_i \qquad\,\,\,\,\, \text{for } |i-j| \geq 2 \\
X_iX_{i+1}X_i &=& X_{i+1}X_iX_{i+1}  \\
X_ie_j &=& e_jX_i \ \qquad\,\,\,\,\, \text{for } |i-j| \geq 2 \\
e_ie_j &=& e_je_i \,\,\,\,\,\,\,\,\, \qquad \text{for } |i-j| \geq 2 \\
X_ie_i &=& e_iX_i \,\,\,=\,\,\, \l e_i  \\
e_iX_{i\pm 1}e_i &=& \l^{-1} e_i \\
X_i X_{i\pm 1} e_i &=& e_{i \pm 1} e_i \,\,\,=\,\,\, e_{i\pm 1} X_i X_{i\pm 1} \\
e_i e_{i\pm 1} e_i &=& e_i \\
%cyclotomic relations
Y^k &=& \sum_{i=0}^{k-1}q_iY^i \\
X_1YX_1Y &=& YX_1YX_1  \\
YX_i &=& X_iY \,\,\,\,\,\,\, \qquad \text{for } i > 1 \\
YX_1Ye_1 &=& \lambda^{-1}e_1 \,\,\,=\,\,\, e_1YX_1Y \\
e_1Y^me_1 &=& A_m e_1 \,\,\,\,\, \qquad \text{for } 0 \leq m \leq k-1.
\end{eqnarray*}
\end{defn}

If the $k^{\rm th}$ order relation $Y^k=\sum_{i=0}^{k-1}q_iY^i$ is omitted, one obtains the affine BMW algebra, studied by Goodman and Hauschild in \cite{GH04}. Thus the cyclotomic BMW algebra is a quotient of the affine BMW algebra. 

Note that the original definition in H\"aring-Oldenburg \cite{HO01} requires that the $k^{\rm th}$ order relation splits over $R$; that is, it can be written 
$\prod_{i=0}^{k-1}(Y-p_i)=0$ for some units $p_i\in R$. We do not require this assumption.

In this paper we shall focus our attention on the algebra $\bt(R)$. We simplify our notation by omitting the index $1$ of $X_1$ and $e_1$. Specifically, $\bt(R)$ is the unital associative $R$-algebra generated by $Y^{\pm 1}, X^{\pm 1}$ and~$e$ 
subject to the following relations: 
\begin{eqnarray}
Y^k &=& \sum_{i=0}^{k-1}q_iY^i \\
X - X^{-1} &=& \d(1-e) \label{eqn:xinv} \\
XYXY &=& YXYX \label{eqn:xyxy} \\
Xe &=& \l e \,\,\,=\,\,\, eX \label{eqn:xele}\\
YXYe &=& \lambda^{-1}e \,\,\,=\,\,\, eYXY \label{eqn:eyxy} \\
eY^me &=& A_m e \,\,\,\,\,\, \qquad \text{for } 0 \leq m \leq 
k-1. \label{eqn:eye}
\end{eqnarray}
Throughout this paper, we make the following assumption:

\begin{assum}\label{NAZD}The element $\delta\in R$ is not a zero divisor.\end{assum}

If we define $q_k = -1$, then $\sum_{l=0}^k q_l Y^l = 0$ and the inverse of $Y$ may
then be expressed as a linear combination of nonnegative powers of $Y$ as follows:
\[
Y^{-1} = - q_0^{-1} \sum_{i=0}^{k-1} q_{i+1}  Y^i.
\]

\noi \textbf{Remark:}  \\
There is an anti-involution $^*$ of $\bt(R)$ such that
\begin{equation} \label{eqn:star}
  Y^* = Y, \quad X^* = X \quad \text{and} \quad e^* = e.
\end{equation}
\\ 
Let $\mathfrak{h}_2^k(R)$ denote the unital associative $R$-algebra with generators $T_0^{\pm 1}$, $T_1^{\pm 1}$ and relations
\begin{eqnarray*}
  T_0 T_1 T_0 T_1 &=& T_1 T_0 T_1 T_0 \\
  T_0^k &=& \sum_{i=0}^{k-1}q_iT_0^i \\
  T_1^2 &=& \d T_1 + 1. 
\end{eqnarray*}
Here $\mathfrak{h}_2^k(R)$ is essentially the Ariki-Koike algebra \cite{AK94}, with the parameters modified. It is clearly a quotient of $\bt(R)$, and indeed $\bt/I \cong \mathfrak{h}_2^k$ as $R$-algebras, where $I$ is the two-sided ideal generated by $e$ in $\bt (R)$.

Using the relations of the algebra and the relation
\begin{eqnarray*}
X Y^l e \= \iy \ix XYX Y^l e \\
&\rel{eqn:xyxy}& \iy \ix Y^l XYX e \\ 
&\stackrel{(\ref{eqn:xele}),(\ref{eqn:eyxy})}{=}& \iy \ix Y^{l-1} e \\
&\stackrel{(\ref{eqn:xinv}),(\ref{eqn:eye})}{=}& Y^{-1} X Y^{l-1} e - \d Y^{l-2} e + \d A_{l-1} \iy e,
\end{eqnarray*}
it is straightforward to show the left ideal of $\bt$ generated by $e$ is $\la{Y^i e \mid 0 \leq i \leq k-1}$, where $\la{S}$ denotes the $R$-submodule spanned by the set $S$.

As a consequence of the results in Goodman and Hauschild \cite{GH04}, the set $\{Y^i e \mid i \in \mathbb{Z}\}$ is linearly independent in the affine BMW algebra, and so it seems natural to expect that the set $\{Y^i e \mid 0 \leq i \leq k-1\}$ will be linearly independent in the cyclotomic BMW algebra. However, for this to hold we must impose additional restrictions on our parameters $A_0, \ldots, A_{k-1}, q_0, \ldots, q_{k-1},q, \l$. We say that our choice of parameters $A_0, \ldots, A_{k-1}, q_0, \ldots, q_{k-1},q, \l$ is \emph{admissible} if these conditions are satisfied. These conditions are given explicitly in Definition \ref{defn:adm} below (under Assumption \ref{NAZD}).

%%%%%%%%%%%%%%%%%%%%%%%%%%%%%%%%%%%%%%%%%%%%%%%%%%%%%%%%%%%%

\section{Construction of $V$}

In this section we construct a $\bt$-module $V$ of rank $k$, under the assumption of admissibility. The action of $\bt$ on $V$ is 
motivated by the action of $\bt$ on the left ideal $\la{Y^i e \mid 0 \leq i \leq k-1}$, and we will show that this left 
ideal is indeed isomorphic to $V$ as $\bt$-modules (see Corollary \ref{BeB}).

Let $V$ be the free $R$-module of rank $k$ with basis $v_0, v_1, \ldots, v_{k-1}$. We want the isomorphism discussed above to send $v_i$ 
to $Y^ie$, so this motivates us to define a linear map $\b{Y}: V \rightarrow V$ by

\begin{eqnarray}
\b{Y}v_i &=& v_{i+1} \text{,\quad for } 0 \leq i \leq k-2, 
\label{eqn:Yv}\\
\b{Y}v_{k-1} &=& \sum_{i=0}^{k-1} q_i v_i.\label{eqn:Yvk-1}
\end{eqnarray}

\noi Since $q_0$ is invertible, this guarantees that $\b{Y}$ is invertible, with inverse
\begin{eqnarray}
\b{Y}^{-1} v_i &=& v_{i-1}\text{,\quad for } 1 \leq i \leq k-1, \label{eqn:iyl}\\
\b{Y}^{-1} v_0 &=& -q_0^{-1} \kum{i} q_{i+1} v_i. \label{eqn:iy0}
\end{eqnarray}
Also $\b{Y}^iv_0=v_i$ for $i = 0, \ldots, k-1$. Define $v_s=\b{Y}^sv_0$ for all integers $s$. The definition of $\b{Y}v_{k-1}$ gives
\[
  \sum_{l=0}^k q_l\b{Y}^lv_0=0.
\]
For any integer $i$, applying $\b{Y}^i$ to this gives
\begin{equation}\label{kthV}
  \sum_{l=0}^k q_l\b{Y}^l v_i = 0.
\end{equation}
Furthermore, since $\b{Y}^l$ is invertible for any integer $l$ and, as $\{v_0, v_1, \ldots, v_{k-1}\}$ is a basis for $V$, the set
\begin{equation}\label{eqn:vbasis}
  \{\b{Y}^lv_0, \b{Y}^lv_1, \ldots, \b{Y}^lv_{k-1}\}
    =\{v_l, v_{l+1}, \ldots, v_{l+k-1}\}
\end{equation}
is also a basis for $V$.

\noi Now let us define linear maps $\b{X}, \b{E}: V \rightarrow V$ by  
\begin{eqnarray}
\b{X}v_0 &=& \l v_0, \label{eqn:Xv0} \\
\b{X}v_1 &=& \l^{-1} \b{Y}^{-1} v_0, \label{eqn:Xv1}\\
\b{X}v_i &=& \b{Y}^{-1} \b{X}v_{i-1} - \d v_{i-2} + \d A_{i-1} \b{Y}^{-1}  v_0 \text{,\quad for } 2 \leq i \leq k-1, \qquad \label{eqn:Xvi} \\
\b{E}v_i &=& A_i v_0, \text{\quad for } 0 \leq i\leq k-1. \label{eqn:Ev}
\end{eqnarray}
Observe that substituting $i=1$ into (\ref{eqn:Xvi}) reproduces (\ref{eqn:Xv1}), as $\l^{-1} = \l - \d + \d A_0$ in $R$. Also, note that the image of the map $\b{E}$ is contained in $\la{v_0}$.
Furthermore, let us denote by $\b{W}$ the map $\b{X} - \d + 
\d\b{E}$. Thus (\ref{eqn:Xv1}) and (\ref{eqn:Xvi}) become 
\begin{equation} \label{eqn:Xvii} 
\b{X} v_i = \b{Y}^{-1} \b{W} 
v_{i-1}\text{,\quad for } 1\leq i \leq k-1. 
\end{equation}

Our aim is to show that $V$ is actually a $\bt$-module, where the action of the 
generators $Y$, $X$ and $e$ are given by the maps $\b{Y}$, 
$\b{X}$ and $\b{E}$, respectively. In order to prove this, we require $\b{Y}^{-1} \b{W}\b{Y}^{-1} = \b{X}$. By equations (\ref{eqn:iyl}) and (\ref{eqn:Xvii}), this relation automatically holds on $v_1, \ldots, v_{k-1}$, so we need only ensure that 
\[(\b{X}-\b{Y}^{-1}\b{W}\b{Y}^{-1})v_0=0.\]

For convenience, we write the left hand side relative to the basis $\{v_0,v_{-1},\ldots,v_{1-k}\}$. Let
\[
  (\b{X}-\b{Y}^{-1}\b{W}\b{Y}^{-1})v_0
    =q_0^{-1}\beta v_0+\sum_{l=1}^{k-1}q_0^{-1}h_lv_{-l},
\]
where $\beta$ and $h_l$ are elements of $R$. Let us calculate $h_l$ and $\beta$ explicitly. By definition, 
since $\l^{-1} = \l - \d + \d A_0$, we have $\b{W}v_0=\l^{-1}v_0$ and
\[
 \b{W}v_i = (\b{X}-\d+\d\b{E})v_i \rel{eqn:Xvii} \b{Y}^{-1}\b{W}v_{i-1}-\d v_i+\d A_iv_0,
\]
for $1\leq i\leq k-1$. 
It is then easy to verify the following by induction on $l$. \\
\textbf{Claim:} \quad For $0 \leq l \leq k-1$, 
\begin{equation} \label{eqn:Wvl}
\b{W}v_l=\l^{-1}v_{-l}+\d\sum_{i=1}^l(A_{l+1-i}v_{1-i}-v_{l-2i+2}).
\end{equation}
Indeed, when $l = 0$, the RHS of equation (\ref{eqn:Wvl}) is simply $\l^{-1} v_0$. Moreover, for $1\leq l\leq k-1$, if the formula holds for $l-1$ then 
\begin{eqnarray*}
 \b{W}v_l
  &=&\b{Y}^{-1}\b{W}v_{l-1}-\d v_l+\d A_lv_0\\
  &\stackrel{\text{ind. hypo.}}{=}&\b{Y}^{-1}\left(\l^{-1}v_{-l+1}+\d\sum_{i=1}^{l-1}(A_{l-i}v_{1-i}-v_{l-2i+1})\right)
      -\d v_l+\d A_lv_0\\
	&=&\l^{-1}v_{-l}+\d\sum_{i=1}^{l-1}(A_{l-i}v_{-i}-v_{l-2i})
      -\d v_l+\d A_lv_0\\
    &=&\l^{-1}v_{-l}+\d\sum_{i=2}^l(A_{l+1-i}v_{1-i}-v_{l-2i+2})
      +\d(A_lv_0-v_l)\\
    &=&\l^{-1}v_{-l}+\d\sum_{i=1}^l(A_{l+1-i}v_{1-i}-v_{l-2i+2}),
\end{eqnarray*}
as required. Thus
\begin{align*}
  -q_0\b{Y}^{-1}\b{W}\b{Y}^{-1}v_0
    &\rel{eqn:iy0} \b{Y}^{-1}\b{W}\sum_{l=1}^kq_lv_{l-1}\\
    &\rel{eqn:Wvl} \b{Y}^{-1}\sum_{l=1}^kq_l
      \left[\l^{-1}v_{-l+1}+\d\sum_{i=1}^{l-1}(A_{l-i}v_{1-i}-v_{l-2i+1})\right]\\
%     &= \sum_{l=1}^kq_l \left[\l^{-1}v_{-l}+\d\sum_{i=1}^{l-1}(A_{l-i}v_{-i}-v_{l-2i})\right]\\
    &\eq{kthV}\l^{-1}\sum_{l=1}^{k-1}q_lv_{-l}-\l^{-1}v_{-k}
      +\d\sum_{l=1}^k\sum_{i=1}^{l-1}q_l(A_{l-i}v_{-i}-v_{l-2i})\\
%     &\rel{kthV}\l^{-1}\sum_{l=1}^{k-1}q_lv_{-l}+\l^{-1}q_0^{-1}\sum_{l=1}^kq_lv_{l-k}
%       +\d\sum_{i=1}^{k-1}\sum_{l=i+1}^kq_l(A_{l-i}v_{-i}-v_{l-2i})\\
    &\rel{kthV}\l^{-1}\sum_{l=1}^{k-1}q_lv_{-l}+\l^{-1}q_0^{-1}\sum_{l=0}^{k-1}q_{k-l}v_{-l}
      +\d\sum_{i=1}^{k-1}\left(\sum_{l=i+1}^kq_lA_{l-i}\right)v_{-i}\\
    &\hphantom{\rel{kthV}} -\, \d\sum_{i=1}^{k-1}\sum_{l=i+1}^kq_lv_{l-2i} \\
    &\eq{kthV}\l^{-1}\sum_{l=1}^{k-1}q_lv_{-l}+\l^{-1}q_0^{-1}\sum_{l=0}^{k-1}q_{k-l}v_{-l}
      +\d\sum_{l=1}^{k-1}\left(\sum_{r=1}^{k-l} q_{r+l}A_r \right)v_{-l}\\
    &\hphantom{\rel{kthV}} -\, \d\sum_{i=1}^{k-1}\sum_{l=i+1}^kq_lv_{l-2i}.
\end{align*}

\noi For now we focus on the last term. Let $k = 2z - \e$, where $\e = 0$ if $k$ is even and $1$ if $k$ is odd. That is, 
let $z :=\lceil k/2\rceil$. Then 
\begin{align*}
  -\d\sum_{i=1}^{k-1}\sum_{l=i+1}^kq_lv_{l-2i}
    &\eq{kthV} -\d \sum_{i=z}^{k-1}\sum_{l=i+1}^kq_lv_{l-2i}
      -\d \sum_{i=1}^{z-1}\sum_{l=i+1}^kq_lv_{l-2i}\\
    &\rel{kthV} -\d \sum_{i=z}^{k-1}\sum_{l=i+1}^kq_lv_{l-2i}
      + \d \sum_{i=1}^{z-1}\sum_{l=0}^iq_lv_{l-2i}\\
    &\eq{kthV} -\d \sum_{i=z}^{k-1}\sum_{l=2i-k}^{i-1}q_{2i-l}v_{-l}
      +\d \sum_{i=1}^{z-1}\sum_{l=i}^{2i}q_{2i-l}v_{-l}\\
    &\eq{kthV} - \d\sum_{l=\e}^{k-2}\sum_{i=\max(l+1,z)}^{\lfloor\frac{l+k}{2}\rfloor} q_{2i-l}v_{-l}      + \d \sum_{l=1}^{k+\e-2}\sum_{i=\lceil\frac{l}{2}\rceil}^{\min(l,z-1)}q_{2i-l}v_{-l}.
\end{align*}

\noi Substituting this into the above expression for $-q_0\b{Y}^{-1}\b{W}\b{Y}^{-1}v_0$, we obtain the following:
\begin{align*}
  \beta v_0+\sum_{l=1}^{k-1}h_lv_{-l}\hspace{-10mm}&\hspace{10mm}=q_0(\b{X}-\b{Y}^{-1}\b{W}\b{Y}^{-1})v_0\\
    &\rel{eqn:Xv0} q_0\l v_0+\l^{-1}\sum_{l=1}^{k-1}(q_l+q_0^{-1}q_{k-l})v_{-l}-\l^{-1}q_0^{-1}v_0\\
    &\hphantom{\rel{eqn:Xv0}}+\d\sum_{l=1}^{k-1}\left(\sum_{r=1}^{k-l}q_{r+l}A_r\right)v_{-l}\\
    &\hphantom{\rel{eqn:Xv0}}-\d\sum_{l=\e}^{k-2}\left(\sum_{i=\max(l+1,z)}^{\lfloor\frac{l+k}{2}\rfloor}\hspace{-5mm}q_{2i-l}\right)v_{-l}
      +\d\sum_{l=1}^{k+\e-2}\left(\sum_{i=\lceil\frac{l}{2}\rceil}^{\min(l,z-1)}\hspace{-5mm}q_{2i-l}\right)v_{-l}.
\end{align*}

\noi Observe that the second last inner sum above is zero when $l=k-1$, as is the last, provided $\e=0$. The upper index of the outer sum can therefore be changed to $k-1$ in both. Now, equating coefficients in the above equation implies that
\begin{eqnarray*}
  \beta&=&q_0\l-q_0^{-1}\l^{-1}+(1-\e)\d,\\
  h_l&=&\l^{-1}(q_l+q_0^{-1}q_{k-l})
    +\d\left[\sum_{r=1}^{k-l}q_{r+l}A_r
    -\sum_{i=\max(l+1,z)}^{\lfloor\frac{l+k}{2}\rfloor}\hspace{-5mm}q_{2i-l}
    +\sum_{i=\lceil\frac{l}{2}\rceil}^{\min(l,z-1)}\hspace{-5mm}q_{2i-l}\right].
\end{eqnarray*}
\begin{defn} \label{defn:adm}
Let $R$ be as in Definition \ref{defn:bnk} (with Assumption \ref{NAZD}). The family of parameters $\left( A_0, \ldots, A_{k-1}, q_0, \ldots, q_{k-1}, q, \l\right)$ is called \textbf{\emph{admissible}} if the following two properties are satisfied:
\begin{enumerate}
\item[(a)] $\l - \l^{-1} = \d(1-A_0)$;
\item[(b)] $h_1=\ldots=h_{k-1}=\beta=0$,
where the $h_l$ and $\beta$ are as defined above.
\end{enumerate}
\end{defn}
\textbf{Remark:} In the case when $\d$ is invertible, the relations $h_l=0$ may be viewed as a triangular system of equations for $A_0,\ldots, A_{k-1}$.

\textbf{Remark:} Note that admissibility is defined in various ways in the literature. The definition of H\"aring-Oldenburg \cite{HO01} differs slightly but involves the same equations as the above. Goodman and Hauschild Mosley \cite{GH109,GH209} use the above definition, under the assumption that $R$ is an integral domain. Rui and Xu \cite{Rui} give a definition involving different equations, which requires additional elements to not be zero divisors but is otherwise equivalent \cite{G10}. See \cite{bnk1} for a more detailed comparison. Note that both these definitions require Assumption \ref{NAZD} to hold. We will show in Corollary \ref{generic} that, under Assumption \ref{NAZD}, the above notion of admissibility is equivalent to $\bt$ being a free $R$-module of rank $3k^2$.

\textbf{Remark:} In \cite{bnk1} we remove Assumption \ref{NAZD}, and modify the above equations to obtain conditions which are necessary and sufficent for $\bt$ to be free of rank $3k^2$; in fact they also ensure that $\bnk$ is free of rank $k^n(2n-1)!!$ for all $n$.
%\pagebreak

Henceforth assume that $R$ is as in Definition \ref{defn:bnk} and the parameters $A_0, \ldots, A_{k-1}, q_0, \ldots, q_{k-1}, q, \l$ are admissible. Denote $\bt(R)$ simply by $\bt$. Now we proceed to show that $V$ is a $\bt$-module. 
 
By our argument above, $\b{Y}^{-1}\b{W}\b{Y}^{-1}=\b{X}$ holds on all of $V$ and hence rearranging gives
\[
\b{Y}\b{X}\b{Y} = \b{W} = \b{X} - \d + \d\b{E}.
\]
Thus $\b{Y}\b{X}\b{Y}\b{X} - \b{X}\b{Y}\b{X}\b{Y} = 
[\b{Y}\b{X}\b{Y},\b{X}] = [\d \b{E}, \b{X}]$, where 
$[\ \text{,}\ ]$ denotes the standard commutator of two maps.

Since $\im(\b{E})\subseteq\la{v_0}$ and $\b{X}v_0=\l v_0$ by definition,
\begin{equation} \label{D1}
\im([\d \b{E},\b{X}]) \subseteq \la{v_0}.
\end{equation}

\noi Let $N = \b{Y}\b{X}\b{Y}\b{X} - 1$. Then 
\be
[N,\b{Y}] &= N\b{Y} - \b{Y}N \\ 
&= \b{Y}\b{X}\b{Y}\b{X}\b{Y} - \b{Y} - \b{Y}\b{Y}\b{X}\b{Y}\b{X} + \b{Y} \\
&= \b{Y}(\b{X}\b{Y}\b{X}\b{Y} - \b{Y}\b{X}\b{Y}\b{X}) \\
&= -\b{Y}[\d \b{E},\b{X}]
\end{align*}
and
\[
[N,\b{Y}^{-1}] = -\b{Y}^{-1} [N,\b{Y}] \b{Y}^{-1} = [\d \b{E},\b{X}] \b{Y}^{-1}.
\]
Therefore, by (\ref{D1}),
\begin{equation} \label{D2}
  \im([N,\b{Y}]) \subseteq \la{v_1} \text{\quad and \quad}  \im([N,\b{Y}^{-1}]) \subseteq \la{v_0}.
\end{equation}
Observe that 
\[
  \b{Y}\b{X}\b{Y}\b{X}v_0 = \l \b{Y}\b{X}\b{Y}v_0 = \l \b{Y}\b{X}v_1 = \l \b{Y} (\l^{-1} \b{Y}^{-1} v_0) = v_0
\]
and
\[
  \b{Y}\b{X}\b{Y}\b{X}v_1 = \b{Y}\b{X}\b{Y} (\l^{-1} \b{Y}^{-1} v_0) = \l^{-1} \b{Y}\b{X} v_0 = \l^{-1} \b{Y} (\l v_0) = v_1,
\]
hence $Nv_0=Nv_1=0$.

\begin{lemma} \label{lem:N}
\begin{equation} \label{D3}
Nv_l \in \la{v_1, v_2, \ldots, v_{l-1}} \tf l \geq 1,
\end{equation}
and
\begin{equation} \label{D4}
Nv_{-m} \in \la{v_0,v_{-1}, \ldots, v_{-(m-1)}} \tf m \geq 0.
\end{equation}
\end{lemma}

\begin{proof}
  We have already established $Nv_0 = Nv_1 = 0$ above. To prove the first assertion, we argue by induction on $l$. Assume that $l \geq 2$ and $Nv_{l-1} \in \la{v_1, v_2, \ldots, v_{l-2}}$. Then, by (\ref{D2}) and the inductive hypothesis,
\[
Nv_l = [N,\b{Y}] v_{l-1} + \b{Y}N v_{l-1} 
\in \la{v_1} + \la{v_2, \ldots, v_{l-1}}.
\]
Thus $Nv_l \in 
\la{v_1, \ldots, v_{l-1}}$ for all $l \geq 1$. \\
Similarly, we prove the second assertion by induction on $m$. Assume $m \geq 1$ 
and $Nv_{-(m-1)} \in \la{v_0,v_{-1}, \ldots, v_{-(m-2)}}$. Then
\be
Nv_{-m} &= [N,\b{Y}^{-1}] v_{-(m-1)} \,+\, \b{Y}^{-1} N v_{-(m-1)} \\
  &\in \la{v_0} + \la{v_{-1}, v_{-2}, \ldots v_{-(m-1)}},
\end{align*}
by (\ref{D2}) and the inductive hypothesis.
So $Nv_{-m} \in \la{v_0,v_{-1}, \ldots, v_{-(m-1)}}$ for all $m \geq 0$.

\end{proof}

\begin{lemma} \label{lem:YXYX}
$\b{Y}\b{X}\b{Y}\b{X}$ is the identity map on $V$.
\end{lemma}
\begin{proof}
We proceed by induction on $i$ to show that $N v_i = 0$ for $i = 0,1,\ldots, k-1$. The cases $i=0$ and $i=1$ have been established above. Suppose that $2 \leq i \leq k-1$ and 
\begin{equation} \label{D5}
  Nv_0 = Nv_1 = \ldots = N v_{i-1} = 0.
\end{equation}
Since $\{v_{i-1}, v_{i-2}, \ldots, v_{i-k}\}$ is a basis for $V$, by (\ref{eqn:vbasis}), we have that
\[
v_i \in \la{v_{i-1},v_{i-2}, \ldots, v_0,v_{-1}, \ldots, v_{i-k}}.
\]
Together with our inductive hypothesis (\ref{D5}), this implies 
\be
Nv_i &\in N \la{v_0,v_{-1}, \ldots, v_{i-k}} \\
&\subseteq \la{v_0,v_{-1}, \ldots, v_{i-k+1}}\text{, by } (\ref{D4}).
\end{align*}
However (\ref{D3}) states that $Nv_i \in \la{v_1, v_2, \ldots, v_{i-1}}$. 
Again using that $\{v_{i-1},v_{i-2}, \ldots, v_{i-k+1},v_{i-k}\}$ forms a basis for $V$, 
\[
\la{v_1, v_2, \ldots, v_{i-1}} \cap \la{v_0,v_{-1}, \ldots, v_{i-k+1}} = 0.
\]
Therefore $Nv_i = 0$. So by induction on $i$, we have proved that $N v_i = 0$ for all $0 \leq i \leq k-1$, as required.
\end{proof}

\begin{lemma} \label{lem:V} (cf. Lemma 25 of H\"aring-Oldenburg \cite{HO01}) \\
Suppose the choice of parameters in $R$ satisfy the admissiblity conditions. Then the following relations 
hold on $V$: 
\begin{eqnarray}
\sum_{l=0}^{k} q_l \b{Y}^l \= 0 \label{eqn:VYk}\\
\b{X}\b{W} \= \b{W}\b{X} \,\,\,=\,\,\, 1 \label{eqn:Vxw}\\
\b{X} \b{E} \= \l \b{E} \,\,\,=\,\,\, \b{E} \b{X} \label{eqn:Vxe}\\
\b{Y}\b{X}\b{Y}\b{X} \= \b{X}\b{Y}\b{X}\b{Y} \label{eqn:Vyxyx}\\
\b{E} \b{Y}^m \b{E} \= A_m \b{E} \tf  \tw{m}, \label{eqn:eyev}\\
\b{E} \b{Y}\b{X}\b{Y} \= \b{Y}\b{X}\b{Y} \b{E}\,\,\,=\,\,\,\l^{-1} \b{E}. 
\label{eqn:Veyxy}
\end{eqnarray} 

\noi Furthermore $V$ is a $\bt(R)$-module, where the 
actions of $Y$, $X$, $\ix$ and $e$ on $V$ are given by the maps 
$\b{Y}$, $\b{X}$, $\b{W}$ and $\b{E}$, respectively.
\end{lemma}

\begin{proof}
The $k^{\rm th}$ order relation (\ref{eqn:VYk}) is immediate from (\ref{kthV}).
As a consequence of Lemma \ref{lem:YXYX},
\[
  \b{W}\b{X}=\b{Y}\b{X}\b{Y}\b{X}=1.
\]
Clearly, this implies
\[
\b{X}\b{W} = \b{X}\b{Y}\b{X}\b{Y} = \iy (\b{Y}\b{X}\b{Y}\b{X}) \b{Y} = \iy (1) \b{Y} = 1.
\]
Hence we have proved (\ref{eqn:Vxw}) and (\ref{eqn:Vyxyx}). Moreover,
\[
[\d\b{E}, \b{X}] = [\b{W},\b{X}] = 0.
\]
We have assumed that $\delta$ is not a zero divisor in $R$, and $\End_R(V)$ is a free $R$-module, so $\b{X}$ commutes with $\b{E}$.
Thus using (\ref{eqn:Ev}) and (\ref{eqn:Xv0}),
\[
\b{E}\b{X} = \b{X}\b{E} = \l \b{E},
\]
proving (\ref{eqn:Vxe}).
Equation (\ref{eqn:eyev}) follows easily from (\ref{eqn:Yv}) and
(\ref{eqn:Ev}). Furthermore, since $\b{Y}\b{X}\b{Y} = \b{W}$,
(\ref{eqn:Vxw}) and (\ref{eqn:Vxe}) imply (\ref{eqn:Veyxy}).
The last assertion of Lemma \ref{lem:V} now follows immediately.
\end{proof}

%%%%%%%%%%%%%%%%%%%%%%%%%%%%%%%%%%%%%%%%%%%%%%%%%%%%%%%%%%%
\section{Constructing the Regular Representation of $\bt$}

In this section, under the assumptions of admissibility, we shall prove 
that $\bt(R)$ is a free $R$-module of rank $3k^2$ and provide an explicit basis of $\bt(R)$; consequently, we show that $\{Y^i e Y^j \mid i,j = 0,\ldots,k-1\}$ is a basis of the two-sided ideal of $\bt(R)$ generated by $e$. 

Recall that $\{v_i \mid r\leq i<r+k\}$ forms a basis for $V$ for any integer~$r$. 
Therefore, the $R$-module $\Xi_0=V\otimes_RV$ has basis
\[
  \{v_{ij} \mid r \leq i < r+k\text{ and }s \leq j < s+k \}
\]
for any integers $r$ and $s$, where $v_{ij}=v_i\otimes v_j$. In particular, 
$\Xi_0$ is $R$-free of rank $k^2$. The anti-involution $^*$ permits us to make 
$\Xi_0 = V \otimes V$ into a $(\bt,\bt)$-bimodule with left action satisfying
\[
  a(v \otimes v') = av \otimes v'
\]
and right action satisfying
\[
  (v \otimes v')b = v \otimes (b^* v'),
\]
for all $a,b \in \bt$ and $v,v' \in V$. In particular, for each integer $j$
there is a $\bt$-module isomorphism from $V$ to $\la{v_{ij} \mid i \in
\mathbb{Z}}$ taking $v_i$ to $v_{ij}$ for all $i$. We now build the regular representation of $\bt$.

Let $\Xi_1$ and $\Xi_2$ be $R$-modules isomorphic to $\Xi_0$, with isomorphisms 
$\xi_1:\Xi_0\rightarrow\Xi_1$ and $\xi_2:\Xi_0\rightarrow\Xi_2$. Define
\[
  \Xi = \Xi_0 \oplus \Xi_1 \oplus \Xi_2,
\]
so that $\Xi$ is $R$-free of rank $3k^2$. For all integers $i$ and $j$, define
$u_{ij} = \xi_1(v_{ij}) \in \Xi_1$ and $w_{ij} = \xi_2(v_{ij}) \in \Xi_2$. Now 
(\ref{kthV}) gives
\[
  \left(\sum_{l=0}^k q_lv_{i+l}\right)\otimes v_j
    =v_i\otimes\left(\sum_{l=0}^k q_lv_{j+l}\right)
    =0.
\]
Expanding and applying $\xi_1$ and $\xi_2$, we obtain
\begin{equation} \label{uvwkth}
 \begin{array}{rcccccl}
  \sum\limits_{l=0}^k q_l v_{i+l,j} &=&
  \sum\limits_{l=0}^k q_l u_{i+l,j} &=&
  \sum\limits_{l=0}^k q_l w_{i+l,j} &=& 0 \\
  \sum\limits_{l=0}^k q_l v_{i,j+l} &=&
  \sum\limits_{l=0}^k q_l u_{i,j+l} &=&
  \sum\limits_{l=0}^k q_l w_{i,j+l} &=& 0, \\
 \end{array}
\end{equation}
for any integers $i$ and $j$.

In order to give $\Xi$ a $\bt$-module structure, we shall define linear maps in $\End_R(\Xi)$, again denoted by $\b{Y}$, $\b{X}$ and $\b{E}$, and show they satisfy the defining relations of $\bt$. We already have a left $\bt$-module structure on $\Xi_0$, so we define $\b{Y}$, $\b{X}$ and $\b{E}$ to act on $\Xi_0$ as $Y$, $X$ and $e$, respectively. 
%\pagebreak

\noi Thus it follows that
\begin{equation}
  \b{Y}v_{ij}=v_{i+1,j}, \label{Yvij}  
\end{equation}
and
\begin{eqnarray*}  
  \b{X}v_{0j}&=&\l v_{0j},\\
  \b{X}v_{ij}&=& \iy X v_{i-1,j} - \d v_{i-2,j} + \d A_{i-1}
\iy v_{0j}, \quad \text{ for }1\leq i\leq k-1,\\
  \b{E}v_{ij}&=&A_iv_{0j}, \quad \text{ for }0\leq i\leq k-1.
\end{eqnarray*}

The aim is to prove $\Xi$ is isomorphic to the regular representation with $v_{ij}$, $u_{ij}$ and $w_{ij}$ 
corresponding to $Y^ieY^j$, $Y^iXY^j$ and $XY^iXY^j$, respectively. Motivated by this, the following definitions are made for $i,j \in \{0,\ldots,k-1\}$:
\begin{eqnarray}
\b{E} u_{ij} \= v_{00} Y^i X Y^j \label{Euij} \\
\b{E} w_{ij} \= v_{00} X Y^i X Y^j\,\,\,=\,\,\,\l \b{E} u_{ij} \label{Ewij} \\
\b{X} u_{ij} \= w_{ij} \label{Xuij} \\
\b{X} w_{ij} \= u_{ij} + \d w_{ij} - \d \l \b{E} u_{ij} \label{Xwij} \\
\b{Y} u_{ij} \= u_{i+1,j} \label{Yuij} \\
\b{Y} w_{ij} \= \b{W} u_{i,j+1} + \d u_{1,i+j} - \d YXv_{ij},
\label{Ywij}
\end{eqnarray}
where, as in the construction of the $\bt$-module $V$, $\b{W}$ denotes the linear map $\b{X} - \d + \d \b{E}$. Note that 
(\ref{Euij}) and (\ref{Ewij}) make use of the right $\bt$-module structure of $\Xi_0$. Since $v_{ij}$, $u_{ij}$ and 
$w_{ij}$ satisfy the $k$-th order relations (\ref{uvwkth}), an easy argument shows that 
(\ref{Yvij})--(\ref{Ywij}) are satisfied for any integers $i$ and $j$.

\begin{lemma}
For any integer $i$,
\begin{equation}\label{eqn:Yw1i}
  \b{Y} w_{1i} = w_{1,i+1}.
\end{equation}
\end{lemma}
\begin{proof}
Because $YXYv_0=X^{-1}v_0=\l^{-1}v_0$, we have
\[
  YXYv_{0i}=\l^{-1}v_{0i}
  \hspace{10mm}\text{and}\hspace{10mm}
  v_{00}YXY=\l^{-1}v_{00}.
\]
Thus
\begin{align*}
  \b{Y}w_{1i}
  &\rel{Ywij} \b{W} u_{1,i+1} + \d u_{1,i+1} - \d YXv_{1i} \\
  &\rel{Yvij} (\b{X}-\d+\d \b{E}) u_{1,i+1} + \d u_{1,i+1} - \d YXY v_{0i} \\
  &\rel{Euij} \b{X} u_{1,i+1} + \d v_{00}YXY^{i+1} - \d \l^{-1} v_{0i} \\
  &\rel{Xuij} w_{1,i+1} + \d \l^{-1} v_{00}Y^i - \d \l^{-1} v_{00}Y^i \\
  &\stackrel{\hphantom{(99)}}{=} w_{1,i+1},
\end{align*}
as required.
\end{proof}

\begin{thm} \label{thm:Xi} 
  The maps $\b{X}$, $\b{W}$, $\b{Y}$, $\b{E}: \Xi \rightarrow \Xi$ defined above satisfy the following identities:
\begin{eqnarray}
  \b{X} \b{W} \= \b{W} \b{X}\,\,\,=\,\,\,1 \label{eqn:XWXi} \\
  \b{X} \b{E} \= \b{E} \b{X}\,\,\,=\,\,\,\l \b{E} \label{eqn:XEXi} \\
  \b{X}\b{Y}\b{X}\b{Y} \= \b{Y}\b{X}\b{Y}\b{X} \label{eqn:XYXYXi}\\
  \b{E}\b{Y}\b{X}\b{Y} \= \b{Y}\b{X}\b{Y}\b{E}\,\,\,=\,\,\,\l^{-1} \b{E}  \label{eqn:eYXYXi} \\
  \b{E} \b{Y}^m \b{E} \= A_m \b{E} \tf \tw{m}, \label{eqn:eYeXi} \\
  \sum_{l=0}^k q_l \b{Y}^l \= 0 \label{Xikth}.
\end{eqnarray}
Furthermore $\Xi$ is $\bt$-module where the actions of the generators $Y$, $X$, $\ix$ and $e$ on $\Xi$ are given by the maps $\b{Y}$, $\b{X}$, $\b{W}$ and $\b{E}$, respectively.
\end{thm}

\begin{proof}
Recall that the actions of $\b{X}$, $\b{Y}$ and $\b{E}$ on $\Xi_0$ are defined by the $\bt$-module structure on $\Xi_0$. 
Thus all the relations hold immediately on $\Xi_0$ and we need only verify them on the $u_{ij}$ and $w_{ij}$.

\bigskip
\noi \textbf{Claim:} \quad $\b{X}\b{W} = \b{W}\b{X} = 1$ and $\b{X} \b{E} = \b{E} \b{X} = \l \b{E}$.

\bigskip
\noi Fix integers $i$ and $j$. It is easy to see that $\b{X}$ and $\b{E}$ preserves the three dimensional $R$-submodule of 
$\Xi$ spanned by $\{u_{ij},w_{ij},\b{E} u_{ij}\}$, and that their actions on this $R$-submodule are given by the matrices
\[
\b{X}' = \left( \begin{array}{ccc}
0 & 1 & 0 \\
1 & \d & 0 \\
0 & -\d\l & \l 
\end{array} \right)
\text{\quad and\quad} 
\b{E}' = \left( \begin{array}{ccc}
0 & 0 & 0 \\
0 & 0 & 0 \\
1 & \l & A_0 
\end{array} \right), 
\]
respectively. Using the fact that $\l - \d + \d A_0 = \l^{-1}$ in $R$, a direct
computation shows
$\b{X}'(\b{X}' - \d I_3 + \d \b{E}')=(\b{X}' - \d I_3 + \d \b{E}')\b{X}'=I_3$.
Thus $\b{X}\b{W}$ and $\b{W}\b{X}$ both act as the identity on this submodule,
and in particular on $u_{ij}$ and $w_{ij}$.  This proves (\ref{eqn:XWXi}). Similarly,
it is easily shown that $\b{X}' \b{E}' = \b{E}' \b{X}' = \l \b{E}'$, so that
$\b{X}\b{E} = \b{E}\b{X} = \l \b{E}$ on $\Xi$, giving~(\ref{eqn:XEXi}).

Now that we have proved (\ref{eqn:XWXi}), we may rewrite $\b{W}$ as $\b{X}^{-1}$. Our next step is to show that relation (\ref{eqn:XYXYXi}) holds on the $u_{ij}$.

\bigskip
\noi \textbf{Claim:} \quad $(\b{X}\b{Y}\b{X}\b{Y})u_{ij} = (\b{Y}\b{X}\b{Y}\b{X})u_{ij}$.

\bigskip
\noi We have
\begin{eqnarray*}
  \b{X}\b{Y}\b{X}\b{Y} u_{ij} &\stackrel{(\ref{Yuij})}{=}& \b{X}\b{Y}\b{X} u_{i+1,j} \\
  &\stackrel{(\ref{Xuij})}{=}& \b{X}\b{Y} w_{i+1,j} \\
  &\stackrel{(\ref{Ywij})}{=}& \b{X}(\b{X}^{-1} u_{i+1,j+1} + \d u_{1,i+j+1} - \d YX v_{i+1,j}) \\
  &\stackrel{(\ref{Xuij}), (\ref{Yvij})}{=}& u_{i+1,j+1} + \d w_{1,i+j+1} - \d XYXY v_{ij}.
\end{eqnarray*}

\noi On the other hand, 
\begin{eqnarray*}
  \b{Y}\b{X}\b{Y}\b{X} u_{ij} 
  &\stackrel{(\ref{Xuij})}{=}& \b{Y}\b{X}\b{Y} w_{ij} \\
  &\stackrel{(\ref{Ywij})}{=}& \b{Y}\b{X} (\b{X}^{-1} u_{i,j+1} + \d u_{1,i+j} - \d YX v_{ij}) \\
  &\rel{Xuij}& \b{Y} (u_{i,j+1} + \d w_{1,i+j} - \d XYX v_{ij}) \\
  &\rel{Yuij}& u_{i+1,j+1} + \d \b{Y} w_{1,i+j} - \d YXYX v_{ij}\\
  &\stackrel{(\ref{eqn:Yw1i}), (\ref{eqn:xyxy})}{=}& u_{i+1,j+1} + \d w_{1,i+j+1} - \d XYXY v_{ij}.
\end{eqnarray*}
This proves our claim:
\begin{equation} \label{eqn:XYXYuij}
(\b{X}\b{Y}\b{X}\b{Y})u_{ij} = (\b{Y}\b{X}\b{Y}\b{X})u_{ij}.
\end{equation}

\noi To prove $(\b{X}\b{Y}\b{X}\b{Y})w_{ij} = (\b{Y}\b{X}\b{Y}\b{X})w_{ij}$, we need to first prove~(\ref{eqn:eYXYXi}).

\bigskip
\noi \textbf{Claim:} \quad $\b{Y}\b{X}\b{Y}\b{E} u_{ij} = \l^{-1} \b{E}u_{ij} = \b{E}\b{Y}\b{X}\b{Y} u_{ij}$ \quad and \quad  $\b{Y}\b{X}\b{Y}\b{E} w_{ij} = \l^{-1} \b{E}w_{ij} = \b{E}\b{Y}\b{X}\b{Y} w_{ij}$.

\bigskip
\noi It is easy to see that 
\begin{align*}
\b{Y}\b{X}\b{Y} \b{E} u_{ij}
 &\rel{Euij} YXYv_{00}Y^iXY^j\\
 &\stackrel{\hphantom{(99)}}{=} \l^{-1}v_{00}Y^iXY^j\\
 &\rel{Euij} \l^{-1} \b{E}u_{ij}.
\end{align*}

\noi Recall that $YXY$ acts as $X-\delta+\delta e$ on $V$. Therefore if $0\leq i\leq k-1$, then
\begin{align*}
\b{E}\b{Y}\b{X}\b{Y} u_{ij}
  &\rel{Yuij} \b{E}\b{Y}\b{X} u_{i+1,j} \\
  &\rel{Xuij} \b{E}\b{Y}w_{i+1,j} \\
  &\rel{Ywij} \b{E}(\b{X}^{-1} u_{i+1,j+1} + \d u_{1,i+j+1} - \d YX v_{i+1,j}) \\
  &\rel{eqn:XEXi} \l^{-1} \b{E}u_{i+1,j+1} + \d \b{E} u_{1,i+j+1} - \d eYX v_{i+1,j} \\
  &\rel{Euij} \l^{-1} v_{00}Y^{i+1}XY^{j+1}+\d v_{00}YXY^{i+j+1}-\d eYX v_{i+1,j}\\
  &\rel{Yvij} \l^{-1} v_{00}Y^i(X-\d+\d e)Y^j+\d\l^{-1}v_{00}Y^{i+j}-\d eYXYv_{ij}\\
  &\kern.18em\stackrel{(\ref{eqn:eyxy})}{=} \l^{-1} v_{00}Y^iXY^j+\d\l^{-1} v_{00}Y^ieY^j-\d\l^{-1} ev_{ij}\\
  &\rel{Euij} \l^{-1}\b{E}u_{ij}+\d\l^{-1} v_{0i}eY^j-\d\l^{-1}A_iv_{0j}\\
  &\stackrel{\hphantom{(99)}}{=} \l^{-1}\b{E}u_{ij}+\d\l^{-1}A_iv_{00}Y^j-\d\l^{-1}A_iv_{00}Y^j\\
  &\stackrel{\hphantom{(99)}}{=} \l^{-1}\b{E}u_{ij}.
\end{align*}

\noi Hence $\b{Y}\b{X}\b{Y}\b{E} u_{ij} = \b{E}\b{Y}\b{X}\b{Y} u_{ij} = \l^{-1} \b{E} u_{ij}$.
Moreover, this implies \[ \b{Y}\b{X}\b{Y}\b{E} w_{ij} \rel{Ewij} \l (\b{Y}\b{X}\b{Y} \b{E} u_{ij})
  = \l \l^{-1} \b{E} u_{ij}  = \l^{-1} \b{E} w_{ij}
\]
and
\begin{align*}
  \b{E}\b{Y}\b{X}\b{Y} w_{ij} &\rel{Xuij} \b{E}\b{Y}\b{X}\b{Y}\b{X} u_{ij} \\
  &\rel{eqn:XYXYuij} \b{E}\b{X}\b{Y}\b{X}\b{Y} u_{ij} \\
  &\rel{eqn:XEXi} \b{X}\b{E} \b{Y}\b{X}\b{Y} u_{ij} \\
  &\stackrel{\hphantom{(19)}}{=}  \l^{-1} \b{X}\b{E} u_{ij} \\
  &\rel{eqn:XEXi} \l^{-1} \b{E}\b{X} u_{ij} \\
  &\rel{Xuij} \l^{-1} \b{E}w_{ij}.
\end{align*}

Hence $\b{Y}\b{X}\b{Y}\b{E} w_{ij} = \l^{-1} \b{E} w_{ij} = \b{E}\b{Y}\b{X}\b{Y} w_{ij}$, finishing our proof of (\ref{eqn:eYXYXi}). We are now ready to show that $\b{X}\b{Y}\b{X}\b{Y} w_{ij} = \b{Y}\b{X}\b{Y}\b{X} w_{ij}$. By (\ref{eqn:XWXi}) and (\ref{eqn:XEXi}), we have $\b{X}^2=1 + \d \b{X} - \d \l \b{E}$. Thus

\begin{eqnarray*}
  \b{Y}\b{X}\b{Y}\b{X} w_{ij} &\rel{Xuij}&\b{Y}\b{X}\b{Y}\b{X}^2 u_{ij}\\
  &=& \b{Y}\b{X}\b{Y} (1 + \d \b{X} - \d \l \b{E}) u_{ij} \\
  &\stackrel{(\ref{eqn:XYXYuij}),(\ref{eqn:eYXYXi})}{=}& (1 + \d \b{X} - \d \l \b{E}) \b{Y}\b{X}\b{Y} u_{ij} \\
&=& \b{X}^2 \b{Y}\b{X}\b{Y} u_{ij} \\
&\rel{eqn:XYXYuij}& \b{X}\b{Y}\b{X}\b{Y}\b{X} u_{ij} \\
&\rel{Xuij}& \b{X}\b{Y}\b{X}\b{Y} w_{ij}.
\end{eqnarray*}
Hence this proves (\ref{eqn:XYXYXi}).

For any integer $l$ and for $\tw{m}$, we have
\[
  \b{E}\b{Y}^mv_{0l}
    =eY^mv_{0l}
    =ev_{ml}
    =A_mv_{0l}.
\]
Thus $\b{E}\b{Y}^m$ acts as $A_m$ on $\la{v_{0l}\mid \tw{l}}$. However, it is clear from the definitions that the image of $\b{E}$ lies in $\la{v_{0l}\mid \tw{l}}$. Therefore $\b{E}\b{Y}^m\b{E}=A_m\b{E}$, 
giving (\ref{eqn:eYeXi}).

It now remains to show that $\b{Y}$ satisfies the $k$-th order relation~(\ref{Xikth}). It follows directly from (\ref{uvwkth}), 
(\ref{Yvij}) and (\ref{Yuij}) that \[\sum_{l=0}^k q_l \b{Y}^l v_{ij} = \sum_{l=0}^k q_l \b{Y}^l u_{ij} =0,\] 
for all integers $i$ and $j$. We shall prove $\sum_{l=0}^k q_l \b{Y}^l w_{ij} =0$ for all $i\geq1$ by induction on $i$.

It is clear from (\ref{eqn:Yw1i}) that $\b{Y}^lw_{1j}=w_{1,j+l}$ for any integer $j$. Thus
\[
  \sum_{l=0}^k q_l \b{Y}^l w_{1j}
    =\sum_{l=0}^k q_l w_{1,j+l}
    \rel{uvwkth}0,
\]
proving the case $i=1$. Suppose the inductive hypothesis holds for some $i \geq 1$.
Recall, by definition (\ref{Xwij}), $\b{X}w_{ij} = u_{ij} + \d w_{ij} - \d \l \b{E} u_{ij}$. Then, using (\ref{Yuij}) and (\ref{Ywij}),
\[
\b{Y}\b{X} w_{ij} = u_{i+1,j} + \d \b{X}^{-1} u_{i,j+1} + \d^2 u_{1,i+j} - \d^2 YX v_{ij} - \d \l \b{Y}\b{E} u_{ij}.
\]
Applying $\b{X}$ to both sides and rearranging gives
\[
  w_{i+1,j} = \b{X}\b{Y}\b{X} w_{ij} - \d u_{i,j+1} - \d^2 w_{1,i+j} + \d^2 XYX v_{ij} + \d \l \b{X}\b{Y}\b{E} u_{ij}.
\]
Now, the terms $u_{i,j+1}$, $XYX v_{ij}$ and $\b{X}\b{Y}\b{E} u_{ij}$ lie in the $R$-submodule $\Xi_0\oplus\Xi_1$. We know the 
operator $\sum_{l=0}^kq_l\b{Y}^l$ acts as the zero operator on this submodule, so applying this operator to both sides gives 
\begin{align*}
  \sum_{l=0}^kq_l\b{Y}^lw_{i+1,j}
    &\stackrel{\hphantom{(41)}}{=} \sum_{l=0}^kq_l\b{Y}^l\b{X}\b{Y}\b{X} w_{ij}-\d^2\sum_{l=0}^kq_l\b{Y}^lw_{1,i+j}\\
    &\rel{eqn:XYXYXi} \b{X}\b{Y}\b{X}\sum_{l=0}^kq_l\b{Y}^l w_{ij}-\d^2\sum_{l=0}^kq_l\b{Y}^lw_{1,i+j}\\
    &\stackrel{\hphantom{(41)}}{=} 0,
\end{align*}
by the inductive hypothesis and the $i=1$ case proved above. Therefore we have proved by induction on $i$ that 
$\sum_{l=0}^k q_l \b{Y}^l w_{ij} =0$, for all $i \geq 1$. However, $\{w_{ij}\mid 1\leq i\leq k\}$ spans $\Xi_2$, so 
(\ref{Xikth}) holds on all of $\Xi$. This completes the proof of the theorem.

\end{proof}

We have now constructed a free $\bt(R)$-module of rank $3k^2$. The following theorem shows that $\Xi$ is in fact isomorphic to $\bt$, naturally considered as a $\bt$-module. Moreover, it states explicitly the bijection between our basis of $\Xi$ and a basis of $\bt$.

\begin{thm} \label{thm:btreg} (cf. Proposition 29 of H\"aring-Oldenburg \cite{HO01}) \\
  The $R$-module homomorphism $\phi: \Xi \rightarrow \bt$ defined by
\begin{eqnarray*}
  v_{ij} &\mapsto& Y^i e Y^j \\
  u_{ij} &\mapsto& Y^i X Y^j \\
  w_{ij} &\mapsto& X Y^i X Y^j
\end{eqnarray*}
is a $\bt$-module isomorphism. Furthermore $\bt$ is $R$-free of rank $3 k^2$ with basis $\{ Y^i e Y^j, Y^i X Y^j, X Y^i X Y^j\mid i,j = 0,\ldots,k-1\}$.
\end{thm}

\begin{proof}
To show that $\phi$ is a $\bt$-module homomorphism, it suffices to check that $\phi(ax_{ij})=a\phi(x_{ij})$ for $a\in\{X,Y,e\}$ and $x\in\{u,v,w\}$. This is straightforward to verify using the algebra relations and the definitions of $\b{X}$, $\b{Y}$ and $\b{E}$.

Let $\k=X^{-1}u_{00}\in\Xi$. Consider the $\bt$-module homomorphism $\zeta: \bt \rightarrow \Xi$ defined by 
\[
a\mapsto a\k.
\]
We will show that $\phi$ and $\zeta$ are mutual inverses. Certainly $\zeta(1)=\kappa$ and
\[
  \phi(\k)=\ix \phi(u_{00})=X^{-1}X=1.
\]
Thus $\phi\zeta(1)=1$. But $\phi\zeta$ is a $\bt$-module homomorphism and $1$ generates $\bt$ as a $\bt$-module, so $\phi\zeta$ is the identity. The same argument will show that $\zeta\phi$ is the identity, provided we can show that $\kappa$ generates $\Xi$, as a $\bt$-module. Let $\bt\k$ be the submodule of $\Xi$ generated by $\k$. Certainly $\bt\k$ contains $X\kappa=u_{00}$. Thus $\bt\k$ also contains $XYu_{00}=w_{10}$. Now, for any integers $i$ and $j$, the element
\[
  Y^{i-1}X^{-1}Y^jw_{10}
    \rel{eqn:Yw1i}Y^{i-1}X^{-1}w_{1j}
    \rel{Xuij}Y^{i-1}u_{1j}
    \rel{Yuij}u_{ij}
\]
is contained in $\bt\k$. Applying $X$, we also have $w_{ij}\in\bt\k$. Finally $\bt\k$ contains
\begin{align*}
  \l^{-1}Y^ieu_{0j}
    &\rel{Euij} \l^{-1}Y^iv_{00}XY^j\\
    &= Y^iv_{00}Y^j\\
    &= v_{ij}
\end{align*}
for all integers $i$ and $j$. Therefore $\bt\k=\Xi$, so that $\phi$ and $\zeta$ are mutual inverses as claimed.
Hence $\bt(R)$ is a free $R$-module of rank $3k^2$ and $\{ Y^i e Y^j, Y^i X Y^j, X Y^i X Y^j\mid i,j = 0, \ldots, k-1 \}$ is an $R$-basis of $\bt(R)$.
  
\end{proof}

\begin{cor} (Lemma 26 of H\"aring-Oldenburg \cite{HO01})\label{BeB} \\
The two-sided ideal of $\bt(R)$ generated by $e$ is a free $R$-module with basis $\{Y^i e Y^j \mid i,j = 0, \ldots, k-1\}$. Moreover, as a left $\bt(R)$-module, it is isomorphic to a direct sum decomposition of $k$ copies of $V$. Specifically
\[
	\la{Y^i eY^j\mid i=0,\ldots,k-1}
\]
is isomorphic to $V$ for each $j=0,\ldots,k-1$.
\end{cor}

\begin{proof}
By Theorem \ref{thm:btreg}, the restriction of $\phi$ to $\Xi_0=V\otimes_RV$ is an injective $\bt$-module homomorphism 
$\phi':V\otimes_RV\hookrightarrow\bt$. Now $V\otimes_RV$ has an $R$-basis $\{v_{ij}\mid i,j = 0, \ldots, k-1\}$, so
\[
  \{\phi'(v_{ij})\mid i,j = 0, \ldots, k-1\}=\{Y^ieY^j\mid i,j = 0, \ldots, k-1\}
\]
is an $R$-basis for $\im(\phi')$. 
Now $\im(\phi')$ is a left ideal of $\bt$, since $\phi'$ is a $\bt$-module homomorphism. But $(Y^ieY^j)^*=Y^jeY^i$. Thus $\im(\phi')$ is invariant under $^*$, so it is a two sided ideal of $\bt$. It clearly contains and is generated by $e$, so it is exactly the two sided ideal generated by~$e$. Finally $V\otimes_RV$ has the following $\bt$-module decomposition:
\[
	V\otimes_RV=\bigoplus_{j=0}^{k-1}V\otimes_Rv_j,
\]
where each summand is isomorphic to $V$ as a $\bt$-module. Applying the isomorphism $\phi'$,
\[
	\bt e\bt=\bigoplus_{j=0}^{k-1}\la{Y^i eY^j\mid i=0,\ldots,k-1},
\]
and again each summand is isomorphic to $V$.
\end{proof}
\begin{cor}\label{generic}
Let $R$ be as in Definition \ref{defn:bnk} (with Assumption \ref{NAZD}). The following are equivalent:
\begin{enumerate}
\item $\bt$ is a free $R$-module of rank $3k^2$.
\item The set $\{Y^ie\mid 0\leq i\leq k-1\}$ is linearly independent in $\bt$.
\item
	$\bt$ admits a module which is free over $R$ of rank $k$, on which the actions of $X$, $Y$ and $e$ are described by 
	(\ref{eqn:Yv}-\ref{eqn:Yvk-1},\ref{eqn:Xv0}-\ref{eqn:Ev}).
\item The parameters of $R$ are admissible.
\end{enumerate}
\end{cor}
\begin{proof}
For any ground ring $R$, the set $\{ Y^i e Y^j, Y^i X Y^j, X Y^i X Y^j\mid i,j = 0,\ldots,k-1\}$ spans $\bt$; indeed the span is 
invariant under left multiplication by $X$, $Y$ and $e$ (the actions are described by (\ref{Euij}-\ref{Ywij})) and contains 
$X^{-1}(Y^0XY^0)=1$. Therefore if (1) holds, this set must be a basis, and in particular (2) holds. If (2) holds, then 
$V=\la{Y^ie\mid 0\leq i\leq k-1}$ is a module satisfying the conditions of (3). We showed that (3) implies (4) in deriving 
the admissibility conditions. Finally (4) implies (1) by Theorem \ref{thm:btreg}.
\end{proof}

%\begin{thm} \label{thm:V}
%The map $\varphi: V \rightarrow \mathrm{span}_R\{Y^i e \mid 0 \leq i \leq k-1\}$ which maps $v_i$ to $Y^i e$ defines a $\bt 
%(R)$-module isomorphism.
%\end{thm}

%\begin{proof}
%Let $U=\mathrm{span}_R\{Y^i e \mid 0 \leq i \leq k-1\}$. It is clear that $\varphi: V \rightarrow U$ is a surjective $R$-module homomorphism. Recall that $U$ is a left ideal in $\bt$. We can therefore define a $\bt$-module homomorphism $\psi: U\rightarrow V$ by
%\[
%  \psi(a)=a(A_0^{-1}v_0), \quad \text{for all } a \in U.
%\]
%Then
%\[
%\psi(\varphi(v_i))=Y^ie(A_0^{-1}v_0)=Y^iv_0=v_i
%\]
%for $0\leq i\leq k-1$, so that $\psi\varphi$ is the identity. Since $\varphi$ is surjective, it follows that $\psi$ and $\varphi$ are inverses. Therefore they are both $\bt$-module isomorphisms.
%\end{proof}

%Hence, as a direct consequence of Theorem \ref{thm:V}, the set $\{Y^i e \mid 0 \leq i \leq k-1\}$ is linearly independent. A more direct proof of this is given as follows. Suppose $\sum_{i=0}^{k-1}\beta_i Y^i e = 0$, where $\beta_i \in R$. Considering the action of both sides on $v_0$ gives
%\[
%  A_0 \sum_{i=0}^{k-1}\beta_i v_i = 0.
%\]
%But $A_0$ is invertible and $v_0, \ldots, v_{k-1}$ are linearly independent over $R$, so each $\beta_i$ must be $0$.

%\bibliography{cyclotomic_BMW_n2_ref}

\bibliographystyle{plain}

\end{document}